\documentclass[finalversion]{pmsart}[2006/01/13]    

\volume{27}%
\fasc{1}%
\years{2007}%
\pages{xxx--xxx}%
\firstpage{1} %
\received{20.6.2006}%
\lastrevision{}

\title{An invariance principle for weakly dependent stationary general models}
\titlethanks{} %
\dedicated{}%
\ShortTitle{Invariance principle for weakly dependent stationary models} %
\ShortAuthors{ P. Doukhan and O. Wintenberger}%
\NumberOfAuthors{2}%

\FirstAuthor{Paul Doukhan}%
\FirstAuthorCity{Malakoff}%
\FirstAuthorAffiliation{Laboratoire de Statistique, CREST, France}%
\FirstAuthorAddress{Timbre J340,\\ 3, avenue Pierre Larousse,\\ 92240 Malakoff, France}%
\FirstAuthorEmail{doukhan@ensae.fr}%
\FirstAutorThanks{}%

\SecondAuthor{Olivier Wintenberger} %
\SecondAuthorCity{Paris} %
\SecondAuthorAffiliation{AMOS, Statistique Appliqu\'ee et
Mod\'elisation Stochastique\\ Universit\'e Paris 1}%
\SecondAuthorAddress{ Centre Pierre
      Mend\`es-France\\ 90 rue de Tolbiac,\\ F-75634 Paris Cedex 13, France}%
\SecondAuthorEmail{olivier.wintenberger@univ-paris1.fr}%
\SecondAutorThanks{}%

\ThirdAuthor{}%
\ThirdAuthorCity{}%
\ThirdAuthorAffiliation{}%
\ThirdAuthorAddress{}%
\ThirdAuthorEmail{}
\ThirdAutorThanks{}%

\FourthAuthor{}%
\FourthAuthorCity{}%
\FourthAuthorAffiliation{}%
\FourthAuthorAddress{}%
\FourthAuthorEmail{}%
\FourthAutorThanks{}%

\FifthAuthor{} %
\FifthAuthorCity{} %
\FifthAuthorAffiliation{}%
\FifthAuthorAddress{} %
\FifthAuthorEmail{}%
\FifthAutorThanks{}%


\MSCcodes{Primary: 60F17.}

\usepackage{color}
\usepackage{xspace}
\def\1{\ensuremath{\mathrm{1}\hspace{-.35em} \mathrm{1}}} 

\def\C{\mathbb{C}}

\def\E{\mathbb{E}}

\def\L{\mathbb{L}}
\def\N{\mathbb{N}}
\def\P{\mathbb{P}}
\def\R{\mathbb{R}}
\def\Z{\mathbb{Z}}
\def\v{\mathop{\rm Var}\nolimits}
\def\cov{\mathop{\rm Cov}\nolimits}
\def\Lip{\mathop{\rm Lip}\nolimits}
\newcommand{\cvg}[1][]{\mathop{\xrightarrow{\text{\upshape\tiny #1}}}}
\renewcommand{\le}{\ensuremath{\leqslant}}
\renewcommand{\ge}{\ensuremath{\geqslant}}
\renewcommand{\epsilon}{\varepsilon}
\newcommand{\ds}{\displaystyle}
\newcommand{\RHS}{right-hand side\xspace}
\newcommand{\refeqn}[1]{eqn$.$~(\ref{#1})\xspace}
\def\resp{resp$.$\xspace}

\def\ie{{\it i.e{$.$}}\xspace}
\def\eg{{\it e.g$.$}\xspace}
\def\iid{independent and identically distributed\xspace}
\newcommand{\Fbar}{\ensuremath{\overline{F}}}
\newcommand{\Gbar}{\ensuremath{\overline{G}}}
\newcommand{\Qbar}{\ensuremath{\overline{Q}}}
\newcommand{\Rbar}{\ensuremath{\overline{R}}}
\newcommand{\xbar}{\ensuremath{\overline{x}}}
\newcommand{\Xbar}{\ensuremath{\overline{X}}}
\newcommand{\ybar}{\ensuremath{\overline{y}}}
\newcommand{\Ybar}{\ensuremath{\overline{Y}}}
\newcommand{\introo}[2]{{\left]{#1,\,#2\,}\right[\kern1pt}}
\newcommand{\introf}[2]{{\left]{#1,\,#2\,}\right]}}
\newcommand{\intrfo}[2]{{\left[{#1,\,#2}\right[\kern1pt}}

\begin{document}
\maketitle
\begin{abstract}
The aim of this article is to refine a weak invariance principle for
stationary sequences given by Doukhan \& Louhichi (1999). Since our conditions are not causal our assumptions need to be stronger than the mixing and causal $\theta$-weak dependence assumptions used in Dedecker \& Doukhan (2003).
Here, if moments of order $>2$ exist, a weak invariance principle and convergence rates in the CLT are obtained;
Doukhan \& Louhichi (1999) assumed  the existence of moments with order $>4$.
Besides the previously used $\eta$- and $\kappa$-weak dependence conditions, we introduce a weaker one, $\lambda$,  which fits the Bernoulli shifts with dependent inputs.
\end{abstract}

%
%
%
\section{Introduction}
Let $(X_t)_{t\in\Z}$ be a real-valued stationary process.
A huge amount of applications make use of such times series.

Several ways of modeling weak dependence have already been proposed.
One of the most popular is the notion of mixing, see \cite{D}
for bibliography; this notion leads to a very nice asymptotic theory, in particular a weak invariance principle under very sharp conditions
(see \cite{Ri} for the strong mixing case).
Such mixing conditions entail restrictions on the model.
For example, Andrews exhibits in \cite{and} the simple counter-example of an
auto-regressive process which does not satisfy any mixing condition
and innovations need much regularity in both MA$(\infty)$ and Markov
models.
Doukhan \& Louhichi introduced in \cite{dl} new weak dependence conditions in order to solve those problems.
We intend to sharpen their assumptions leading to a weak invariance principle.

A common approach to derive a weak invariance principle for stationary sequences is based on a martingale difference approximation. This approach was first explored by Gordin in \cite{gordin};
necessary and sufficient conditions were found by Heyde in \cite{heyde}.
Let ${\cal M}_t$ be a filtration. Heyde's martingale difference approximation is equivalent to the existence of moments of order 2 and
\begin{equation}
\label{condheyde}
\sum_{t=0}^\infty\E(X_t|{\cal M}_0)\quad\mbox{ and
}\quad\sum_{t=0}^\infty(X_{-t}-\E(X_{-t}|{\cal M}_0))\qquad\mbox{
converge in }\L^2.
\end{equation}%
%
%
Martingale theory leads directly to invariance principles (see
also \cite{VOLN93}).
In the following, the adapted case refers to the special case where $X_t$ is ${\cal M}_t$-measurable.
The natural filtration is written as
${\cal M}_t=\sigma(Y_i,i\le t)$
for \iid inputs $(Y)_{t\in\Z}$; thus $X_t$ can be written as a function of the past inputs:
\begin{equation}%
\label{bercaus}
X_t=H(Y_t,Y_{t-1},\ldots).
\end{equation}
Then only the first series in (\ref{condheyde}) needs to be considered. Using the Lindeberg technique, Dedecker \& Rio relax (\ref{condheyde}) in \cite{dedrio}.
Bernstein's blocks method allowed Peligrad \& Utev to also improve on (\ref{condheyde}) in \cite{PU2}.
Such projective conditions are related to dependence coefficients; Dedecker \& Doukhan obtain sharp results for the causal $\theta$-dependence in \cite{ded} and Merlev\`{e}de {\it et al.} address the  mixing cases in a nice survey paper \cite{MP}.

Martingale difference approximation is not always easy, for instance in the particular case where a natural filtration does not exist.
The most striking example is given by associated sequences $(X_t)_{t\in\Z}$. Let us recall this notion. A series is said to be associated if $\cov(f_1,f_2)\ge 0$ for any two coordinatewise nondecreasing functions $f_1$ and $f_2$ of $(X_{t_1},\dots,$ $X_{t_m})$
with $\v(f_1)+\v(f_2)<\infty$.
However, Newman \& Wright obtain in \cite{NW} a weak invariance principle under the existence of second order moments and
\begin{equation}
\sigma^2=\sum_{t\in\Z}\cov(X_0,X_t)<\infty.
\end{equation}
Theorems \ref{th3} and \ref{th2} propose invariance principles
under general assumptions:
they apply to the non-causal Bernoulli shifts with weakly dependent inputs
$(Y_t)_{t\in\Z}$,
\begin{equation}%
\label{bernoncaus}%
X_t=H(Y_{t-j},j\in\Z).
\end{equation}

Heredity of weak dependence through such non-linear functionals follows from a new $\lambda$-weak dependence property;
a function of a $\lambda$-weak dependence process is $\lambda$-weakly
dependent, see Section \ref{bernsection}.
Analogous models with dependent inputs are already considered by \cite{BBD}. If $X_t=\sum_{j\in\Z}\alpha_jY_{t-j}$,  Peligrad \& Utev prove in \cite{PU1} that the Donsker invariance principle holds for $X$ as soon as it holds for the innovation process $Y$.
The non-linearity of $H$ considered here is an important feature which has not been frequently discussed in the past.
The condition of moments with order $>4$ on the observations  needed in \cite{dl} is reduced to a one of moments with order $>2$ and the results rely on specific decays of the dependence coefficients.
We do not reach the second order moment condition of \cite{heyde} (or projective
conditions) and \cite{NW}.
We \emph{conjecture}  that some times series satisfy  weak dependence conditions with fast enough decay rates in order to ensure a Donsker type theorem but they do not satisfy neither condition (\ref{condheyde}) nor other projective criterion (see \cite{MP}) nor association nor Gaussianity.
More general models (\ref{bernoncaus}) are considered here while causal models (\ref{bercaus}) fit to the adapted case and to projective conditions.
However, proving this conjecture is really difficult since condition (\ref{condheyde}) has to be checked for each $\sigma$-algebra ${\cal M}_0$.

The paper is organized as follows.
In Section \ref{sectionresult} we introduce various weak dependent coefficients in order to state our main results.
Section \ref{sectionexample} is devoted to examples of weak dependent models for which we discuss our results.
We shall focus on examples of $\lambda$-weakly dependent sequences. Proofs are given in the last section; we first derive conditions ensuring the convergence of the series $\sigma^2$.
A bound of the $\Delta$-moment of a sum (with $2<\Delta<m$) is proved in Section~\ref{smom}; this bound is of an independent interest since {\it eg.} it directly yields the strong laws of large number.
The standard Lindeberg method with Bernstein's blocks is developed in Section~\ref{sectiontheo12} and yields our versions of the Donsker theorem. Convergence rates of the CLT are obtained in Section~\ref{rate}.
\section{Definitions and main results}\label{sectionresult}
\subsection{Weak dependence assumptions}
\begin{definition} [Doukhan \& Louhichi, 1999] The process $(X_t)_{t\in\Z}$ is said to
be $(\epsilon,\psi)$-weakly dependent if there exist a sequence
$\epsilon(r)\downarrow0$ (as $r\uparrow\infty$) and a function
$\psi:\N^2\times(\R^+)^2\to\R^+$ such that
$$
\left|\cov(f(X_{s_1},\ldots,X_{s_u}),g(X_{t_1},\ldots,X_{t_v})\right|\le
\psi(u,v,\Lip f,\Lip g)\epsilon(r),
$$
for any $r\ge 0$ and any $(u+v)$-tuples such that $s_1\le\cdots\le
s_u\le s_u+r\le t_1\le\cdots\le t_v$, where the real valued
functions $f,g$ are defined respectively on $\R^u$ and $\R^v$,
satisfy $\|f\|_\infty\le1$, $\|g\|_\infty\le1$ and are such that $\Lip f+\Lip\,
g<\infty$ where
$$
\Lip f=\sup_{(x_1,\ldots,x_u)\ne(y_1,\ldots,y_u)}
\frac{\left|f(x_1,\ldots,x_u)-f(y_1,\ldots,y_u)\right|}{|x_1-y_1|+\cdots+|x_u-y_u|}
$$
\end{definition}
Specific functions $\psi$ yield notions of weak dependence
appropriate to describe various examples of models:
\begin{itemize}
\item $\kappa$-weak dependence for which $\psi(u,v,a,b)=uvab$,
in this case we simply denote $\epsilon(r)$ as\ $\kappa(r)$;
\item $\kappa'$ (causal)  weak dependence for which $\psi(u,v,a,b)=vab$,
in this case we simply denote $\epsilon(r)$ as $\kappa'(r)$;
this is the causal counterpart of $\kappa $ coefficients which is recalled only for completeness;
\item $\eta$-weak dependence, $\psi(u,v,a,b)=ua+vb$,
in this case we write $\epsilon(r)=\eta(r)$ for short;
\item $\theta$-weak dependence is a causal dependence which refers to the function $\psi(u,v,a,b)$ $=vb$,
in this case we simply denote $\epsilon(r)=\theta(r)$ (see \cite{ded});
this is the causal counterpart of $\eta $ coefficients which is recalled only for completeness;
\item $\lambda$-weak dependence $\psi(u,v,a,b)=uvab+ua+vb$,
in this case we write $\epsilon(r)=\lambda(r)$.
\end{itemize}
\begin{commentaryRemark}
Besides the fact that it includes $\eta$ and $\kappa$-weak
dependences,
this new notion of $\lambda$-weak dependence will be
proved to be convenient,
for example, for the Bernoulli shifts with associated inputs (see Lemma \ref{th3} below).
\end{commentaryRemark}
\begin{commentaryRemark}
If  functions $f$ and $g$ are complex-valued,
the previous inequalities remain true if we substitute $\epsilon(r)/2$ to $\epsilon(r)$.
A useful case  of such complex-valued functions is
$f(x_1,\ldots,x_u)=\exp\left(it(x_1+\cdots+x_u)\right)$
for each
$t\in\R,$ $u\in\N^*$ and $(x_1,\ldots,x_u)\in\R^u$ (see Section~\ref{sectiontheo12}).
This indeed corresponds to the characteristic function adapted to derive the convergence in distribution.
\end{commentaryRemark}
\subsection{Main results}
Let $(X_t)_{t\in\Z}$ be a real-valued stationary sequence of mean 0 satisfying
\begin{equation}%
\label{mo2+}
\E |X_0|^m<\infty, \qquad \mbox{ for a real number }m>2.
\end{equation}
Let us assume that
\begin{eqnarray} \label{sigma2}
\sigma^2=\sum_{k\in\Z}\cov(X_0,X_k)=\sum_{k\in\Z}\E X_0 X_k\ge 0.
\end{eqnarray}
Denote by $W$ the standard Brownian motion and by $W_n$ the partial sums process
\begin{equation}%
\label{prodonsker}%
W_n(t)=\frac1{\sqrt n}\sum_{i=1}^{[nt]}X_i, \quad
\mbox{for } t\in[0,1],\ n\ge 1.
\end{equation}
We now present our main results,
which are new versions of the Donsker weak invariance principle.
\begin{theorem}[$\kappa$-dependence]%
\label{cltkappa}%
\label{th1}
Assume that the 0-mean $\kappa$-weakly dependent stationary process
$(X_t)_{t\in\Z}$ satisfies \refeqn{mo2+} and $\kappa(r)={\cal
O}(r^{-\kappa})$ (as $r\uparrow\infty$) for $\kappa>2+\frac1{m-2}$
then the previous expression $\sigma^2$ is finite and
$$
W_n(t)\cvg[D]_{n\to\infty}\sigma W(t),\quad\mbox{in the Skorohod space } D([0,1]).
$$
\end{theorem}
\begin{commentaryRemark}
Under the more restrictive $\kappa'$ condition, Bulinski \& Sashkin obtain
invariance principles with the sharper assumption
$\kappa'>1+\frac1{m-2}$ in \cite{Bul}. Our loss is explained by the fact that
$\kappa'$-weakly dependent sequences satisfy $\kappa'(r)\ge
\sum_{s\ge r}\kappa_s$. This simple bound directly follows from
the definitions. \end{commentaryRemark}

The following result relaxes the previous dependence assumptions at the price of a faster decay of the dependence coefficients.
\begin{theorem}%
[$\lambda$-dependence]%
\label{cltepsilon}%
\label{th2}
Assume that the 0-mean $\lambda$-weakly dependent stationary inputs
satisfies \refeqn{mo2+} and  $\lambda(r)={\cal O}(r^{-\lambda})$ (as
$r\uparrow\infty$) for $\lambda>4+\frac2{m-2}$ then $\sigma^2$
is finite and
$$
W_n(t)\cvg[D]_{n\to\infty}\sigma W(t),\quad\mbox{in the Skorohod space } D([0,1]).
$$
\end{theorem}
\begin{commentaryRemark}
We do not achieve better results for $\eta$ or $\theta$-weak dependence cases than the one for $\lambda-$dependence.
In comparison with the result obtained by \cite{ded},
our results are not as good under $\theta$-weak dependence.
We work under more restrictive moment conditions than these authors.
The same remark applies for all projective measures of dependence;
here, we refer to results by \cite{heyde}, \cite{NW}, \cite{dedrio} and \cite{PU2}.
\end{commentaryRemark}
\begin{commentaryRemark}
However, the  example of Section~\ref{bernsection} stresses the fact that such results are not systematically better than those of Theorem \ref{th2};
for such general examples, we even conjecture that theorems of \cite{heyde}, \cite{NW}, \cite{dedrio} or \cite{PU2} do not apply.
\end{commentaryRemark}
\begin{commentaryRemark}
The technique of the proofs is based on the Lindeberg method.
In fact, we prove that
$ \left|\E \left(\phi(S_n/\sqrt n)-\phi(\sigma
N)\right)\right|=o\left(n^{-c}\right)$ ($\phi$
denotes here the characteristic function) for $0<c<c^*$
where $c^*$ depends only on the parameters $m$ and $\kappa$ or $\lambda$ respectively.
If $m$ and $\kappa$ (or $\lambda$) both tend to infinity,
we notice that $c^*\to\frac14$. As $\kappa$ or $\lambda$ tends to infinity and $m<3$, $c^*$ always remains smaller than $(m-2)/(2m-2)$ (see Proposition \ref{prop3} in Section~\ref{rate} for more details).
\end{commentaryRemark}
\begin{commentaryRemark}
Using a smoothing lemma also yields an analogous bound for the
uniform distance
$$
\sup_{x\in\R}\left|\P\left(\frac1{\sqrt n}S_n\le
x\right)-\P\left(\sigma N\le x\right)\right|=
o\left(n^{-c}\right),\qquad \mbox{ for some }c<c'.
$$
A first and easy way to control $c'$ is to let $c'=c^*/4$ but the
corresponding rate is a really bad one (see \eg in \cite{dp}).
The Esséen inequality holds with the optimal exponent $\frac12$ in the \iid case (see \cite{Petrov}) and \cite{Ri} reaches the exponent $\frac13$ in the case of strongly mixing sequences. In Proposition \ref{prop3} of Section~\ref{rate}, we achieve
$c'>c^*/4$. Analogous results have been
settled in
\cite{dlang} for weakly dependent random fields. Previous results by \cite{heyde}, \cite{NW}, \cite{dedrio} or \cite{PU2} do not derive such convergence rates for the Kolmogorov distance.
\end{commentaryRemark}
Let us denote by
$\R^{(\Z)}=\bigcup_{I>0}\left\{z\in\R^\Z\big/~z_i=0,\ |i|>I\right\}$, the set of finite sequences of real numbers.
We consider functions $H:\R^{(\Z)}\to \R$ such that if
$x,y\in\R^{(\Z)}$ coincide for all indexes but one, let say $s\in\Z$, then
\begin{equation}
\label{as}
 |H(x)-H(y)|\le b_s(\|z\|^\ell\vee 1)|x_s-y_s|
\end{equation}
where $z\in\R^{(\Z)}$ is defined by $z_s=0$ and $z_i=x_i=y_i$ if
$i\ne s$.
Here $\|x\|=\sup_{i\in\Z}|x_i|$.
In Section \ref{bernsection}, we prove the existence of 
$X_n=\lim_{I\to\infty}H\big((Y_{n-j}\1_{\{j\le
I\}})_{j\in\Z}\big)$
where $(Y_t)_{t\in\Z}$ is a weakly dependent real-valued input process.
We denote this process by $X_n=H(Y_{n-j},j\in\Z)$ for simplicity
and we derive its $\lambda$-weak invariance properties.
Various asymptotic results follow,
among which our weak invariance principle, Theorem \ref{th2}.
\begin{corollary}%
\label{th3}%
Let $(Y_t)_{t\in\Z}$ be a stationary $\lambda$-weakly dependent
process (with dependence coefficients $\lambda_{Y}(r)$) and
$H:\R^{(\Z)}\to \R$ satisfying the condition given by (\ref{as}) for some $\ell\ge 0$. Let us assume that there exist real numbers $m,\ m'$ with
$\E|Y_0|^{{m'}}<\infty$ such that $m>2$ and
$m'\ge(\ell+1)m$.

Then $X_n=H(Y_{n-i},{i\in\Z})$ exists and satisfies the weak invariance principle in the following cases:
\begin{itemize}
\item \textbf{Geometric case:} $b_r\le Ce^{-b|r|}$ and
$\lambda_{Y}(r)\le De^{-ar}$ for $a,b,C,D>0$.
\item \textbf{Riemannian case:} If $b_r\le C(1+|r|)^{-b}$ for some
$b>2$ and $\lambda_{Y}(r)\le Dr^{-a}$ for $a,C,D>0$ with
\begin{eqnarray}
\label{donsklin} a&>&\frac{1+b}{b-1}\left(4+\frac{2}{m-2}\right),\ \mbox{ if }\ \ell=0,\ b>1;\\
\nonumber
a&>&\frac{b(m'-1+\ell)}{(b-2)(m'-1-\ell)}\left(4+\frac{2}{m-2}\right),\
\mbox{ if }\ \ell>0,\ b>2.
\end{eqnarray}
\end{itemize}
\end{corollary}
\begin{commentaryRemark}
The previous conditions are also tractable in the mixed cases. We
explicitly state them for $\ell>0$:
\begin{itemize}
\item $b_r\le Ce^{-b|r|}$, $\lambda_{Y}(r)\le Dr^{-a}$, if moreover  $\displaystyle
a>\frac{m'-1+\ell}{m'-1-\ell}\left(4+\frac{2}{m-2}\right)$ and $b,C,D>0$.
\item $b_r\le C|r|^{-b}$ and $\lambda_{Y}(r)\le De^{-ar}$, for $a,C,D>0$ with $\displaystyle
b>\frac{6m-10}{m-2}.$
\end{itemize}
\end{commentaryRemark}
\section{Examples}\label{sectionexample}
Theorem \ref{th3} is useful to derive the weak invariance principle
in various cases.
This section is aimed at a detailed treatment of the Bernoulli shifts with dependent inputs.
The important class of Lipschitz functions of dependent inputs is presented in a separate section. The importance of our results is highlighted by the models of the first subsection.
More general non-linear models are considered in the second subsection.
Some of those examples illustrate the conjecture we made in the introduction but we were not able to formally prove it.
\subsection{Lipschitz processes with dependent inputs}%
\label{seclin}
Consider Lipschitz functions $H:\R^{(\Z)}\to\R$, \ie such
that \refeqn{as} is satisfied for $\ell=0$.
A simple example of this situation is the two-sided linear sequence
\begin{equation}
\label{linear} X_t=\sum_{i\in\Z}\alpha_iY_{t-i}
\end{equation}
with dependent inputs $(Y_t)_{t\in\Z}$. As quoted by \cite{LASO00} for the case of linear processes with dependent input there exists a very general solution;
essentially any Donsker type theorem for the stationary inputs implies the central limit theorem for any linear process driven by such inputs.
More precisely, Theorem 5 of \cite{PU1} states that this process even satisfies the Donsker invariance principle if
$\sum_{j}|\alpha_j|<\infty$.

A simple example of Lipschitz non-linear functional of dependent
inputs is
\begin{equation}
\label{linearabs}
X_t=\left|\sum_{i\in\Z}\alpha_iY_{t-i}\right|-\E\left|\sum_{i\in\Z}\alpha_iY_{-i}\right|
\end{equation}
In this case the inequality (\ref{as}) holds with $\ell=0$ and
$b_r\le|\alpha_r|$.

 Another example of this situation is the
following stationary process
$$X_t=Y_t\left(a+\sum_{j\ne0}a_jY_{t-j}\right)-\E
Y_t\left(a+\sum_{j\ne0}a_jY_{t-j}\right), $$
where the inputs
$(Y_t)_{t\in\Z}$ are bounded. In this case, the inequality (\ref{as})
also holds with $\ell=0$ and $b_s\le 2\|Y_0\|_\infty |a_s|$.

To apply our result, we compute the weak dependence coefficients of such models.
\begin{lemma}\label{lemlipsch}
Let $(Y_t)_{t\in\Z}$ be a strictly stationary  process with a finite moment of order $m\ge 1$ and $H:\R^{(\Z)}\to \R$ satisfying the condition (\ref{as}) for $\ell= 0$ and some nonnegative sequence $(b_s)_{s\in\Z}$
such that $L=\sum_jb_j<\infty$.
Then,
\begin{itemize}
\item
the process
$X_n=H(Y_{n-j},j\in\Z):=\lim_{I\to\infty}H\big(Y_{n-j}\1_{\{j\le I\}},j\in\Z\big)$ is a strictly stationary process with finite moments of order $m$.
\item
if the input process $(Y_t)_{t\in\Z}$ is $\lambda$-weakly dependent (the weak dependence coefficients are denoted by $\lambda_{Y}(r)$),
then $(X_t)_{t\in\Z}$ is  $\lambda$-weakly dependent with
$$\lambda(k)= \inf_{2r\le
k}\left[2\sum_{|i|\ge
r}b_i\|Y_0\|_1+(2r+1)^2L^2\lambda_{Y}(k-2r)\right].$$
\item if the input process $(Y_t)_{t\in\Z}$ is $\eta$-weakly dependent (the weak dependence coefficients are denoted by $\eta_{Y}(r)$) then $(X_t)_{t\in\Z}$ is  $\eta$-weakly dependent and
$$
\eta(k)=\inf_{2r \le k} \left[ 2\sum_{|i|\ge
r}b_i\|Y_0\|_1+(2r+1)L\,\eta_{Y}(k-2r)\right].
$$
\end{itemize}
\end{lemma}
\begin{commentaryRemark}
Let $(Y_t)_{t\in\Z}$ be a strictly stationary  process with a
finite moment of order $m>2$. If $L=\sum_j|\alpha_j|<\infty$, the
process $X_n=\sum_{j\in\Z}\alpha_jY_{n-j}$ is a strictly
stationary process with finite moments of order $m$ which
satisfies the assumptions of
Lemma~\ref{lemlipsch} with $b_j=|\alpha_j|$. Even if the weak invariance principle is
already given in \cite{PU1} our result is of an independent
interest, {\it e.g.} for functional estimation purposes. For
non-linear Lipschitz functionals it yields new central limit
theorems. \end{commentaryRemark} The result of Theorems \ref{th1}
and \ref{th2} holds systematically in geometric cases. Then it is
assumed Riemannian decays, \ie there exists $\alpha, C>0 $ such
that
$$b_r\le C r^{-\alpha}.$$The conditions from \cite{heyde} are compared below with the conditions of
Theorems \ref{th1} and \ref{th2} for specific classes of inputs $(Y_t)_{t\in\Z}$.

\subsubsection{LARCH($\infty$) inputs}
A vast literature is devoted to the study of conditionally heteroskedastic models.
A simple equation in terms of a vector-valued process allows a unified treatment of those models,
see \cite{dtw}.
Let  $( \xi_t )_{t \in \Z}$ be an \iid centered real-valued sequence and $a,  a_j , {j \in \N^{\ast}}$ be  real numbers.
LARCH($\infty$) models are solutions of the recurrence equation
\begin{equation}%
\label{archv1}
Y_t = \xi_t \left( a + {\sum_{j = 1}^\infty} a_j Y_{t -j} \right).
\end{equation}
We provide below sufficient conditions for the following chaotic expansion
\begin{equation}%
\label{archv2}
  Y_t = \xi_t \left( a + \sum_{k = 1}^{\infty}
 \sum_{j_1, \ldots, j_k \ge 1}a_{j_1}\xi_{t
  - j_1} a_{j_2}\cdots a_{j_k} \xi_{t
  - j_1 - \cdots - j_k} a\right).
\end{equation}
Assume that  $\Lambda=\| \xi_0 \|_{m}\sum_{j \ge 1} | a_j |< 1$
then one (essentially unique) stationary solution of \refeqn{archv1} in  $\L^{m}$ is given by \refeqn{archv2}. This solution is $\theta$-weakly dependent with $ \theta_Y(r) \le  Kr^{1-a} \log^{a -1}r $ for some constant $K>0$.
This implies the same bound on their coefficients $(\lambda_Y(r))_{r\ge0}$.
Condition (\ref{donsklin}) gives the weak invariant principle for $(X_t)_{t\in\Z}$ under the conditions $ \E| \xi_0 |^m< +\infty$ for $m>2$, $\alpha>1$, and
$$\displaystyle
a>\frac{1+\alpha}{\alpha-1}\left(4+\frac{2}{m-2}\right)+1.$$
The model (\ref{linear}) is also an Heyde's martingale difference approximation (\ref{condheyde}) as soon as
$$\sum_{k\ge1}\sqrt{\sum_{i\ge k}\alpha_i^2}<+\infty.$$
Necessary conditions for weak invariance principle follow as $\alpha>3/2$,
$| a_j |\le Cj^{-a}$ for some $a>1$, $ \E \xi_0^2< +\infty,$ and
$\| \xi_0 \|_2\sum_{j \ge 1} | a_j |< 1$.
These conditions are not optimal since in this case the process is adapted to the filtration ${\cal M}_t=\sigma(\xi_i,i\le t)$.
Peligrad \& Utev extend in \cite{PU2} the Donsker theorem to the cases where $\alpha>1/2$.
Thus, our conditions are not optimal compared to those of \cite{PU1} in the linear case as in \refeqn{linear}. However, for non-linear Lipschitz functional, the result seems to be new.
\subsubsection{Non-causal LARCH($\infty$) inputs}
\label{noncausalarch}
The previous  approach  extends for the case of non-causal LARCH($\infty$) inputs
$$
Y_t = \xi_t \left( a + \sum_{j\ne0} a_j Y_{t -j} \right).
$$
Doukhan {\it et al.} prove in \cite{dtw} the same results of existence as for the previous causal case (just replace summation over $j>0$ by summation over $j\ne0$) and the dependence becomes of the $\eta$ type with
$$\eta(r) =  \left(\| \xi_0 \|_\infty
  \sum_{0\le 2k <r } k \Lambda^{k - 1} A \left( \frac{r}{2k} \right)
  + \frac{\Lambda^{r/2}}{1 - \Lambda} \right)\E | \xi_0 | | a |$$
where
$ A(x)=\sum_{|j|\ge x}|a_j|,$ $ \Lambda=\| \xi_0 \|_\infty\sum_{j \ge 1} | a_j | <1 .$
>From condition (\ref{donsklin}) the weak invariance principle holds for $(X_t)_{t\in\Z}$ if  $\| \xi_0
\|_\infty<\infty$, $\alpha>1$ with
$$\displaystyle
a>\frac{1+\alpha}{\alpha-1}\left(4+\frac{2}{m-2}\right)+1.$$
Notice that a very restrictive new assumption is that inputs need to be uniformly bounded in this non-causal case.
This result is new, a conjecture is that (\ref{condheyde}) does not hold.
\subsubsection{Non-causal, non-linear inputs}%
\label{exampleopt}
The weak dependence properties of non-causal and non-linear inputs $Y_t$ are recalled, see \cite{dl} for more details.
Let $H:{(\R^d)}^\Z\to \R$ be a measurable function.
If the sequence $(\xi_n)_{n\in \Z}$ is independent and identically distributed on $\R^d$, the Bernoulli shift with input process
$(\xi_n)_{n\in \Z}$ is defined as
$$
Y_n=H\left((\xi_{n-i})_{i\in \Z}\right), \qquad n\in \Z.
$$
Such Bernoulli's shifts are $\eta$-weakly dependent (see
\cite{dl}) with $ \eta(r)\le 2\delta_{[r/2]} $ if
\begin{equation}%
\label{regul}
\E\left|H\left(\xi_{j},j\in\Z\right)-
H\left(\xi_{j}\1_{|j|\le r},j\in\Z\right)\right|\le\delta_r.
\end{equation}
Then condition (\ref{donsklin}) leads to the invariance principle for $(X_t)_{t\in\Z}$ if $\E|Y_0|^m<\infty$ for $m>2$, $\alpha>1$ and $\delta_r\le Kr^{-\delta}$ for
$$\delta>\frac{1+\alpha}{\alpha-1}\left(4+\frac{2}{m-2}\right).$$
Conditions (\ref{condheyde}) of \cite{heyde} do not give clear conditions on coefficients for these models. We do not know other weak invariance principle in that general context.

\subsubsection{Associated inputs}

A process is associated  if\ $\cov
\left(f(Y^{(n)}),g(Y^{(n)})\right)\ge0$ for any coordinatewise non-decreasing
function $f,g:\R^n\to \R$ such that the previous covariance makes
sense with $Y^{(n)}=(Y_1,\ldots,Y_n)$. The $\kappa$-weak dependence
condition is known to hold for associated or Gaussian sequences. In both cases this condition
will hold with
$$
\kappa(r)=\sup_{j\ge r}|\cov(Y_0,Y_j)|
$$
Notice the absolute values are needed only in the second case since
for associated processes these covariances are nonnegative. Independent sequences as well are associated and Pitt proves in \cite{Pitt} that a Gaussian process with nonnegative covariances
is also associated. Finally, we recall that non-decreasing functions
of associated sequences remain associated.  Associated models are
classically built this way from \iid sequences, see \cite{L2001}.

Suppose that the inputs $(Y_t)_{t\in\Z}$ are such that $\kappa(r)\le
C r^{-a}$ (for some $a,C>0$). For the associated cases and model (\ref{linear}), the invariance principle of \cite{NW} follows from  remark of \cite{LS}
as soon as $\E Y^2<+\infty$, $a>1$ and $\alpha>1$.
These conditions are optimal, they correspond to $\sum_j\cov(X_0,X_j)<\infty$.
Such strong conditions are due to the fact that zero correlation implies independence for associated
processes.
Our conditions for invariance principle are much stronger: $\E |Y|^m<+\infty$ with $m>2$, $\alpha>1$ and
$$a>\frac{1+\alpha}{\alpha-1}\left(4+\frac{2}{m-2}\right).$$
For non-linear Lipschitz cases as in \refeqn{linearabs} the result seems to be new. In the special case of $\kappa$-weak dependent inputs that are not associated,
the optimal weak invariance principle of \cite{NW} does not apply, see \cite{dl} for examples.
\subsection{The Bernoulli shifts with dependent inputs}\label{bernsection}
Let $H:\R^{(\Z)}\to \R$ be a measurable (non necessarily Lipschitz) function and $X_n=H(Y_{n-i},i\in\Z)$.
Such models are proved to exhibit either $\lambda$ or $\eta$-weak dependence properties.
Because the Bernoulli shifts of $\kappa$-weak dependent inputs are neither $\kappa$ nor $\eta$-weakly dependent,
the $\kappa$ case is here included in the $\lambda$ one.

Consider the non-Lipschitz function $H$ defined by
$$
H(x)=\sum_{k=0}^K\sum_{j_1,\dots,j_k}a^{(k)}_{j_1,\ldots,j_k}x_{j_1}\cdots
x_{j_k}.
$$
In this case, Lemma \ref{lemlipsch} does not apply.
To derive weak dependence properties of such processes,
we assume that $H$ satisfies the condition (\ref{as}) with $\ell\ne0$,
which remains a stronger assumption than for the case of independent inputs,
see \refeqn{regul}. Relaxing Lipschitz assumption on $H$ is possible if we assume the existence of higher moments for the inputs.
The following lemma gives both the existence and the weak dependence properties of such models
\begin{lemma}
\label{lambd}
Let $(Y_i)_{i\in\Z}$ be a stationary process and $H:\R^{(\Z)}\to \R$ satisfies condition (\ref{as}) for some $\ell> 0$ and some sequence $b_j\ge0$ such that $\sum_j|j|b_j<\infty$.
Let us assume that there exist a pair of real numbers $(m,m')$ with $\E|Y_0|^{{m'}}<\infty$ such that $m\ge1$ and $m'\ge (\ell+1)m$.
Then,
\begin{itemize}
\item the process $X_n=H(Y_{n-i},i\in\Z)$ is well defined in $\L^m$: it is a strictly stationary process;
\item if the input process $(Y_i)_{i\in\Z}$ is $\lambda$-weakly dependent
(the weak dependence coefficients are denoted by $\lambda_{Y}(r)$),
then $X_n$ is  $\lambda$-weakly dependent and there exists a
constant $c>0$ such that
$$\lambda(k)=c\inf_{r\le
[k/2]}\left[ \sum_{|j|\ge
r}|j|b_j+(2r+1)^2\lambda_{Y}(k-2r)^{\frac{{m'}-1-\ell}{{m'}-1+\ell}}\right];$$
\item if the input process $(Y_i)_{i\in\Z}$ is $\eta$-weakly dependent
(the weak dependence coefficients are denoted by $\eta_{Y}(r)$)
then $X_n$ is  $\eta$-weakly dependent and there exists a constant $c>0$ such that
$$
\eta(k)=c\inf_{r\le [k/2]} \left[\sum_{|j|\ge r}|j|b_j+ (2r+1)^{
1+\frac{\ell}{{m'}-1} }\eta_{Y}(k-2r)^{\frac{{m'}-1-\ell}{{m'}-1}}\right].
$$
\end{itemize}
\end{lemma}
Such models where already mentioned in the mixing case by \cite{B} and  \cite{BBD}.
The proofs are deferred to Section \ref{prooflemma}.

\subsubsection{Volterra models with dependent inputs}
Consider the function $H$ defined by
$$
H(x)=\sum_{k=0}^K\sum_{j_1,\dots,j_k}a^{(k)}_{j_1,\ldots,j_k}x_{j_1}\cdots
x_{j_k},
$$
then if $x,y$ are as in \refeqn{as}
\begin{multline*}
H(x)-H(y)=\sum_{k=1}^K\sum_{u=1}^k
\sum_{\tiny\begin{array}{c}j_1,\ldots,j_{u-1}\\
j_{u+1},\ldots,j_k
\end{array}}a^{(k)}_{j_1,\ldots,j_{u-1},s,j_{u+1},\ldots,j_k}\times\\
x_{j_1}\cdots x_{j_{u-1}}(x_s-y_s) x_{j_{u+1}}\cdots x_{j_{k}}.
\end{multline*}
From the triangular inequality we thus derive that the previous lemma \ref{lambd} may be written with $\ell=K-1$,
$$b_s=
\sum_{k=1}^K {\sum}^{(k,s)}|a^{(k)}_{j_1,\ldots,j_k}|
$$
where $\sum^{(k,s)}$ stands for the sums over all indices in $\Z^k$
where one of the indices $j_1,\ldots,j_k$ takes on the value $s$ and
$$
L\equiv
\sum_{k=0}^K\sum_{j_1,\ldots,j_k}\left|a^{(k)}_{j_1,\ldots,j_k}\right|.
$$
For example,
$|a^{(k)}_{j_1,\ldots,j_k}|\le C \left(j_1\vee\cdots\vee
j_k\right)^{-\alpha}$
or
$\le C \exp\left(-\alpha(j_1\vee\cdots \vee j_k)\right)$
respectively yield $b_s\le C's^{d-1-\alpha}$ or $b_s\le
C'e^{-\alpha s}$ for some constant $C'>0$.
\subsubsection{Markov stationary inputs} Markov stationary sequences satisfy a recurrence equation
$$
Z_n=F(Z_{n-1},\ldots,Z_{n-d},\xi_n)$$
where $(\xi_n)$
is a sequence of \iid random variables. In this case $Y_n=(Z_n,\ldots,Z_{n-d+1})$ is a Markov chain
$Y_n=M(Y_{n-1},\xi_n)$ with
\begin{equation}
M(x_1,\ldots,x_d,\xi)=(F(x_1,\ldots,x_d,\xi),x_1,\ldots,x_{d-1}).
\label{duflomarkov}
\end{equation}
 Theorem 1.IV.24 of \cite{duflo} proves that \refeqn{duflomarkov} has a stationary solution $(Z_n)_{n\in\Z}$  in $\L^m$ for $m\ge1$ as soon as $\|F(0,\xi)\|_m<\infty$ and there exist a norm $\|\cdot\|$ on $\R^d$ and a real number $a \in \intrfo{0}{1}$ such that  $\|F(x,\xi)-F(y,\xi)\|_m\le a\|x-y\|$.
In this setting $\theta$-dependence holds with $\theta_{Z}(r)={\cal O}(a^{r/d})$ (as $r\uparrow\infty$).
We shall not give more details about the significative examples provided in \cite{dlsv}.
Indeed, we already mentioned that our results are sub-optimal in such causal cases;
such dependent sequences may however also be used as inputs for the Bernoulli shifts.
\subsubsection{Explicit dependence rates}
We now specify the decay rates from Lemma \ref{lambd}.
For standard decays of the previous sequences, it is easy to get the following explicit bounds.
Here $b,c,C,D,\lambda,\eta>0$ are constants which may differ from one case to the other.
\begin{itemize}
\item
If $b_j\le C(|j|+1)^{-b}$  and $\lambda_{Y}(j)\le Dj^{-\lambda}$,
\resp $\eta_{Y}(j)\le Dj^{-\eta}$ then from a simple calculation, we
optimize both terms in order to prove that
$\lambda(k)\le ck^{-\lambda\left(1-\frac2b\right)\frac{{m'}-1-\ell}{{m'}-1+\ell}}$,
\resp $\eta(k)\le ck^{-\eta\frac{(b-2)({m'}-2)}{(b-1)({m'}-1)-\ell}}$.\\
Note that in the case where ${m'}=\infty$ this exponent may be arbitrarily close to $\lambda$ for large values of $b>0$.
This exponent may thus take all possible values between 0 and $\lambda$.
\item If $b_j\le Ce^{-|j|b}$ and $\lambda_{Y}(j)\le De^{-j\lambda}$,
respectively $\eta_{Y}(j)\le De^{-j\eta}$, we obtain $\lambda(k)\le ck^{2}e^{-\lambda
k\frac{b({m'}-1-\ell)}{b({m'}-1+\ell)+2\eta({m'}-1-\ell)}}$, \resp
$\eta(k)\le ck^{ \frac{{m'}-1-\ell}{{m'}-1}}e^{-\eta
k\frac{b({m'}-2)}{b({m'}-1)+2\eta({m'}-2)}}$.\\
The geometric decay of both $(b_j)_{j\in\Z}$ and the weak dependence coefficients of the inputs ensure the geometric decay of the weak dependence coefficients of the Bernoulli shift.
\item If we assume that the coefficients  $(b_j)_{j\in\Z}$ associated with the Bernoulli shift have a geometric decay,
say $b_j\le Ce^{-|j|b}$ and that $\lambda_{Y}(j)\le Dj^{-\lambda}$
(resp{$.$} $\eta_{Y}(j)\le Dj^{-\eta}$)  we obtain the bounds
$\lambda(k)\le ck^{-\lambda \frac{{m'}-1-l}{{m'}-1+\ell}}\log^{2} k $, \resp
$\eta(k)\le ck^{-\eta \frac{{m'}-2}{{m'}-1}}\log^{1+\frac{\ell}{{m'}-1}} k  $.\\
If ${m'}=\infty$ tightness is reduced by a factor $\log^2 k$ with respect to the dependence coefficients of the input dependent series $(Y_t)_{t\in\Z}$.
\item If we assume that the coefficients  $(b_j)_{j\in\Z}$ associated with the Bernoulli shift have a Riemannian decay, say $b_j\le C(|j|+1)^{-b}$ and that
$\lambda_{Y}(j)\le De^{-j\lambda}$ (\resp $\eta_{Y}(j)\le D e^{-j\eta}$)  we find
$\lambda(k)\le ck^{2-b}$, \resp $\eta(k)\le ck^{2-b}$.
\end{itemize}

All models or functions of models we present here are $\lambda$-weakly dependent.
We treat some basic examples in detail when a discussion with other results is possible.
We believe that for some models,
$\lambda$-weak invariance properties follow from easy computations,
and then, statistical results like our weak invariance principle.
\section{Proofs of the main results}\label{proofs}
Our proof for central limit theorems is based on a truncation method.
For a truncation level $T\ge1$ we shall denote $\Xbar_k=f_T(X_k)-\E f_T(X_k)$ with $f_T(X)=X\vee (-T)\wedge T$.
>From now on, we shall use the convenient notation $a_n\preceq b_n$ for two real sequences $(a_n)_{n\in\N}$ and $(b_n)_{n\in\N}$ when there exists some constant $C>0$ such that $|a_n|\le Cb_n$ for each integer $n$.
We also remark that $\Xbar_k$ has moments of all orders because it is bounded.
In the entire sequel, we denote $\mu=\E|X_0|^m$.
For any $a\le m$, we control the moment $\E|f_T(X_0)-X_0|^a$ with Markov inequality
$$\E|f_T(X_0)-X_0|^a\le \E |X_0|^a\1_{\{|X_0|\ge T\}}\le \mu T^{a-m},$$
thus using Jensen inequality yields
\begin{equation}%
\label{tronc}
\|\Xbar_0-X_0\|_a\le 2\mu^{\frac1a}T^{1-\frac ma}.
\end{equation}
Starting from this truncation, we are now able to control the limiting variance as well as the higher order moments.

In this section we prove that the central limit theorems  corresponding to the convergence $W_n(1)\to W(1)$ in both Theorems \ref{th1}
and \ref{th2} hold and we shall provide convergence rates corresponding to these central limit theorems.
The weak invariance principle is obtained in a standard way from such central limit theorems and tightness, which follows from Lemma \ref{lambd}, by using the classical Kolmogorov-Centsov tightness criterion, see \cite{B}.
In the last subsection, we prove Lemma \ref{momlambda} that states the properties of
our (new)  Bernoulli's shifts with dependent inputs.
\subsection{Variances}
\begin{lemma} [Variances]
\label{variance} If one of the following conditions holds\begin{eqnarray}\label{varkappa}
\sum_{k=0}^\infty\kappa(k)\quad &<&\infty
\\
\label{varepsilon}
\sum_{k=0}^\infty\lambda(k)^{\frac{m-2}{m-1}}&<&\infty
\end{eqnarray}
then the series $\sigma^2$ is convergent.
\end{lemma}
\begin{proof} Using the fact that $\Xbar_0=g_T(X_0)$ is a function of
$X_0$ with $\Lip g_T=1$ snd $\|g_T\|_\infty\le 2T$, we derive
\begin{equation}
\label{covtronc}
|\cov (\Xbar_0,\Xbar_k)|\le \kappa(k) \mbox{ or }
(4T+1)\lambda(k), \mbox{ respectively.}
\end{equation}
In the $\kappa$ dependent case, the truncation may thus be omitted
and
\begin{equation}%
\label{covkappa}
|\cov (X_0,X_k)|\le\kappa(k).
\end{equation}
In the following, we shall only consider $\lambda$ dependence. We develop
$$
\cov (X_0,X_k)=\cov (\Xbar_0,\Xbar_k)+\cov
(X_0-\Xbar_0, X_k)+\cov (\Xbar_0,X_k-\Xbar_k).
$$
We use a truncation $T$ (to be determined) and the two previous
bounds \refeqn{tronc} and \refeqn{covtronc}; then the H\"older inequality
with the exponents $1/a+1/m=1$ yields
\begin{eqnarray*}
|\cov (X_0,X_k)|&\le&  (4T+1)\lambda(k)+2\|X_0\|_m\|\Xbar_0-X_0\|_a
\\
&\le&  (4T+1)\lambda(k)+4\mu^{1/a+1/m}T^{1-m/a}
\\
&\le&   (4T+1)\lambda(k)+4\mu T^{2-m}.
\end{eqnarray*}
Chosing $T^{m-1}=\mu/\lambda(k)$ we obtain
\begin{equation}%
\label{covepsilon}
|\cov
(X_0,X_k)|\le
9\mu^{\frac1{m-1}}\lambda(k)^{\frac{m-2}{m-1}}.
\end{equation}
\end{proof}
\subsection{A $\Delta$-order moment bound}\label{smom}
\begin{lemma}\label{momlambda}
Let $(X_t)_{t\in\Z}$ be a stationary and centered process. Let us assume that $\E |X_0|^m<\infty$, and that this process is either $\kappa$-weakly dependent with $\kappa(r)={\cal
O}\left(r^{-\kappa}\right)$ or  $\lambda$-weakly dependent with
$\lambda(r)={\cal O}\left(r^{-\lambda}\right)$. If $\kappa>2+\frac{1}{m-2}$, or $\lambda>4+\frac{2}{m-2},$ then for
all $\Delta>2$ small enough there exists a constant $C>0$ such that
$$
\|S_n\|_\Delta\le C\sqrt{n}.
$$
\end{lemma}
\begin{commentaryRemark}
$\Delta\in\introo{2}{2+A\wedge B\wedge1}$ where $A$ and $B$
are constants smaller than $m-2$ and depend on $m$ and
respectively $\kappa$ or $\lambda$. Equations (\ref{A}) and
(\ref{B}) precise the previous involved constants $A$ and $B$.
\end{commentaryRemark}
\begin{commentaryRemark}
The constant satisfies $\displaystyle
C>\left(\frac{5}{2^{(\Delta-2)/2}-1}\right)^{1/\Delta} \sum_{k\in
\Z} |\cov(X_0,X_k)|.$ Under the conditions of this lemma, Lemma
\ref{variance} entails
$$
 c\equiv\sum_{k\in \Z} |\cov(X_0,X_k)|<\infty.
$$
\end{commentaryRemark}
\begin{commentaryRemark}
The result is sketched from \cite{Bul}. However, their dependence condition is of a causal nature while our is not.
It explains a loss with respect to the exponents $\lambda$ and $\kappa$.
In their $\kappa'$-weak dependence setting the best possible value of the exponent is 1 while it is 2 for our non-causal dependence.
\end{commentaryRemark}
\noindent {\bf Proof of Lemma \ref{momlambda}.} For convenience, let denote in the sequel $\Delta=2+\delta$ and $m=2+\zeta$. Like in \cite{Ibr3} or \cite{Bul}, we proceed by induction on $k$ for $n\le 2^k$ to prove
that
\begin{equation}%
\label{rec}
\|1+|S_n|\|_\Delta\le C\sqrt{n}.
\end{equation}
We assume that \refeqn{rec} is satisfied for all $n\le 2^{K-1}$. Setting $N=2^K$ we have to find a bound for $\|1+|S_N|\|_\Delta$.
We can divide the sum $S_N$ into three blocks:
the first two blocks have the same size $n\le 2^{K-1}$ and are denoted by $Q$ and $R$; the third block $V$, located between $Q$ and $R$, has cardinality $q<n$.
We then have $\|1+|S_N|\|_\Delta\le \|1+|Q|+|R|\|_\Delta+\|V\|_\Delta$.
The term $\|V\|_\Delta$ is directly bounded by $\|1+|V|\|_\Delta\le
C\sqrt{q}$ from the recurrence assumption.
Writing $q=N^b$ with $b<1$,
then this term is of order strictly smaller than $\sqrt N$.
For $\|1+|Q|+|R|\|_\Delta$, we have
\begin{eqnarray*}
\E(1+|Q|+|R|)^\Delta &\le &\E(1+|Q|+|R|)^2(1+|Q|+|R|)^\delta,\\
&\le& \E(1+2|Q|+2|R|+(|Q|+|R|)^2)(1+|Q|+|R|)^\delta.
\end{eqnarray*}
We expand the \RHS of this expression; the following terms appear
\begin{itemize}
\item $\E(1+|Q|+|R|)^\delta\le 1+|Q|_2^\delta+|R|_2^\delta\le 1+2c^\delta(\sqrt n )^\delta$,
\item $\E |Q|(1+|Q|+|R|)^\delta\le \E |Q|((1+|R|)^\delta+|Q|^\delta)$\\
$\le \E |Q|(1+|R|)^\delta+\E|Q|^{1+\delta}$.

The term $\E|Q|^{1+\delta}$ is bounded by $\|Q\|_2^{1+\delta}$ and
then by $c^{1+\delta}(\sqrt n) ^{1+\delta}$.
The term $\E|Q|(1+|R|)^\delta$ is bounded by
$\|Q\|_{1+\delta/2}\|1+|R|\|_\Delta^\delta$ using H\"older inequality. It is at least of
order $cC^\delta(\sqrt{n})^{1+\delta}$, analogous to the latter one, where we exchange the roles of $Q$ and $R$.
\item $\E (|Q|+|R|)^2(1+|Q|+|R|)^\delta$. For this term,
we use an inequality from \cite{Bul}
\begin{multline*}
\E (|Q|+|R|)^2(1+|Q|+|R|)^\delta\\
\le \E|Q|^{\Delta}+\E|R|^{\Delta}
+5(\E Q^2(1+|R|)^\delta+\E R^2(1+|Q|)^\delta).
\end{multline*}
Now
$\E|Q|^{\Delta}\le C^\Delta(\sqrt{n})^\Delta$ is bounded by using
\refeqn{rec}.
The second term is its analogous with $R$ substituted to $Q$.
The third term has to be handled with a particular care, as
follows.
\end{itemize}
We use the weak dependence notion to control $\E Q^2(1+|R|)^\delta$ and $\E R^2(1+|Q|)^\delta$. Denote by $\Xbar$ the variable $X\vee
T\wedge(- T)$ for a real $T>0$ to be determined later. By extension $\Qbar$ and $\Rbar$ denote the truncated sums of the
variables $X_i$. We have
\begin{multline*}
\E |Q|^2(1+|R|)^\delta\le\\
\E Q^2||R|-|\Rbar||^\delta+\E (Q^2-\Qbar^2)(1+|\Rbar|)^\delta
+\E \Qbar^2(1+|\Rbar|)^\delta.
\end{multline*}
We begin with a control of $\E Q^2||R|-|\Rbar||^\delta$.
Using the H\"older inequality with $2/m+1/m'=1$ yields
$$
\E Q^2||R|-|\Rbar||^\delta\le \|Q\|_m^2\|||R|-|\Rbar||^\delta\|_{m'}
$$
$\|Q\|_\Delta$ is bounded using \refeqn{rec} and
$$
||R|-|\Rbar||^{\delta m'}
\le |R|^{\delta m'}\1_{\{|R|>T\}}\le |R|^{\delta
m'}\1_{|R|>T}.
$$
We then bound $\1_{|R|>T}\le (|R|/T)^\alpha$ with
$\alpha=m-\delta m'$, hence
$$\E||R|-|\Rbar||^{\delta m'}\le E|R|^m T^{\delta m'-m}.
$$
By convexity and stationarity, we have $\E|R|^m\le n^m\E|X_0|^m$,
so that
$$
\E Q^2(|R|-|\Rbar|)^\delta\preceq n^{2+m/m'}T^{\delta-m/m'}.
$$
Finally, remarking that $m/m'=m-2$, we obtain
$$
\E Q^2(|R|-|\Rbar|)^\delta\preceq n^{m}T^{\Delta-m}.
$$
We get the same bound for the second term
$$
\E (Q^2-\Qbar^2)(1+|\Rbar|)^\delta\preceq n^{m}T^{\Delta-m}.
$$
For the third one, we introduce a covariance term
$$
\E \Qbar^2(1+|\Rbar|)^\delta\le \cov(\Qbar^2,(1+|\Rbar|)^\delta)
+\E\Qbar^2\E(1+|\Rbar|)^\delta.
$$
The latter is bounded with $|Q|_2^2|R|_2^\delta\le c^\Delta\sqrt{n}^\Delta$.
The covariance is controlled as follows by using weak-dependence
\begin{itemize}
\item in the $\kappa$-dependent case: $n^2T\kappa(q)$,
\item in the $\lambda$-dependent case: $n^3T^2\lambda(q)$.
\end{itemize}
We then choose either the truncation
$T^{m-\delta-1}=n^{m-2}/\kappa(q)$ or
$T^{m-\delta}=n^{m-3}/\lambda(q)$.
At this point, the three terms of the decomposition are of the same order
\begin{eqnarray*}
\E |Q|^2(1+|R|)^\delta&\preceq &\left(n^{3m-2\Delta}\kappa(q)^{m-\Delta}\right)^{1/(m-\delta-1)},\ \mbox{under }\kappa\text{-dependence},\\
\E |Q|^2(1+|R|)^\delta&\preceq&
\left(n^{5m-3\Delta}\lambda(q)^{m-\Delta}\right)^{1/(m-\delta)}, \
\mbox{under }\lambda\text{-dependence}.
\end{eqnarray*}
Let $q=N^b$, we note that $n\le N/2$ and this term is of order
$N^{\frac{3m-2\Delta+b\kappa(\Delta-m)}{m-\delta-1}}$ under
$\kappa$-weak dependence
 and the order $N^{\frac{5m-3\Delta+b\lambda(\Delta-m)}{m-\delta}}$ under $\lambda$-weak dependence.
Those terms are thus negligible with respect to $N^{\Delta/2}$ if
\begin{eqnarray}
\label{hypkap}&&\kappa>\frac{3m-2\Delta-\Delta/2(m-\delta-1)}{b(m-\Delta)},\ \mbox{under }\kappa\text{-dependence},\\
\label{hyplam}&&\lambda>\frac{5m-3\Delta-\Delta/2(m-\delta)}{b(m-\Delta)},\
\mbox{under }\lambda\text{-dependence}.
\end{eqnarray}
Finally, using this assumption, $b<1$ and $n\le N/2$, we derive the
bound for some suitable constants $a_1,a_2>0$
$$
\E(1+|S_N|)^\Delta\le \left(2^{-\delta/2}C^\Delta+5\cdot2^{-\delta/2}c^\Delta+a_1N^{-a_2}\right)  \left(\sqrt{N}\right)^\Delta.
$$
Using the relation between $C$ and $c$, we conclude that \refeqn{rec}
is also true at the step $N$ if the constant $C$ satisfies
$2^{-\delta/2}C^\Delta+5\cdot2^{-\delta/2}c^\Delta+a_1N^{-a_2}\le
C^\Delta$. Choose
$C>\left(\frac{5c^\Delta+a_12^{\delta/2}}{2^{\delta/2}-1}\right)^{1/\Delta}$
with $c=\sum_{k\in \Z} |\cov(X_0,X_k)|$, then the previous relation
holds.
 Finally, we use \refeqn{hypkap} and \refeqn{hyplam}
to find a condition on $\delta$. \\
In the case of $\kappa$-weak dependence, we rewrite inequality
\refeqn{hypkap} as
$$
0>\delta^2+\delta(2\kappa-3-\zeta)-\kappa\zeta+2\zeta+1.
$$
It leads to the following condition on $\delta$
\begin{equation}%
\label{A}
\delta<\frac{\sqrt{(2\kappa-3-\zeta)^2+4(\kappa\zeta-2\zeta-1)}+\zeta+3-2\kappa}{2}\wedge 1=A.
\end{equation}
We do the same in the case of the $\lambda$-weak dependence
\begin{equation}%
\label{B}
\delta<\frac{\sqrt{(2\lambda-6-\zeta)^2+4(\lambda\zeta-4\zeta-2)}+\zeta+6-2\lambda}{2}\wedge 1=B. \quad \square
\end{equation}
\begin{commentaryRemark} The bounds $A$ and $B$ are always smaller than
$\zeta$. \end{commentaryRemark}
\subsection{Proofs of Theorems \ref{th1} and \ref{th2}}
\label{sectiontheo12}
Let $S=\frac1{\sqrt n}S_n$ and consider $p=p(n)$ and $q=q(n)$ in
such a way that
$$
\lim_{n\to\infty} \frac1{q(n)}=\lim_{n\to\infty} \frac
{q(n)}{p(n)}=\lim_{n\to\infty} \frac {p(n)}n=0
$$
and $k=k(n)=n/[p(n)+q(n)]$
$$ Z=\frac1{\sqrt{n}}\left(U_1+\cdots+U_k\right),\qquad \mbox{ with }
U_j=\sum_{i\in B_j}X_i$$
where
$
B_j=\introf{(p+q)(j-1)}{(p+q)(j-1)+p}\cap\N
$
is a subset of $p$ successive integers from $\{1,\dots,n\}$ such
that, for $j\neq j'$, $B_j$ and $B_{j'}$ are at least distant of
$q=q(n)$ from each other. We denote by $B'_j$ the block between $B_j$ and $B_{j+1}$ and
$V_j=\sum_{i\in B_j'}X_i$. $V_k$ is the last block of $X_i$ between
the end of $B_k$ and $n$.
Furthermore, let
$\sigma_p^2=\v(U_1)/p=\sum_{|i|<p}(1-|i|/p)\E X_0X_i$,
and let
$$
Y=\frac{U'_1+\dots+U'_k}{\sqrt{n}}, \qquad U'_j\sim {\cal N}(0,p\,\sigma_p^2)
$$
where the Gaussian variables $V_j$ are mutually independent and also independent of the sequence $(X_n)_{n\in\Z}$.
We also consider a sequence $U_1^*,\ldots,U_k^*$ of mutually independent random variables with the same distribution as $U_1$ and we let
$Z^*=\left(U_1^*+\cdots+U_k^*\right)/{\sqrt{n}}$.
In the entire section, we fix $t\in\R$ and we define $f:\R\to\C$ by
$f(x)=\exp\left\{it x\right\}$. Then $\displaystyle \E f(S)-f(\sigma
N)=\E f(S)-f(Z)+\E f(Z)-f(Z^*)+\E f(Z^*)-f(Y)+\E f(Y)-f(\sigma N). $
Lindeberg method is devoted to prove that this expression converges
to 0 as $n\to\infty$. The first and the last terms in this inequality
are referred to as the auxiliary terms in this Bernstein-Lindeberg
method. They come from the replacement of the individual initial  --
non-Gaussian and Gaussian respectively -- random variables by their
block counterparts. The second term is analogous to that obtained
with decoupling and turns the proof of the central limit theorem to
the independent case. The third term is referred to as the main term
and following the proof under independence it will be bounded above
by using a Taylor expansion. Because of the dependence structure, in
the corresponding bounds, some additional covariance terms will
appear.

\noindent The following subsections are organized as follows:
we first consider the auxiliary terms and the main terms are then
decomposed by the usual Lindeberg method
and the corresponding terms coming from the dependence or the usual remainder terms (standard for the independent case) are considered in separate subsections.
In the last one, we collect these calculations to obtain the central limit theorem.
\subsubsection{Auxiliary terms}%
\label{auxil}
Using Taylor expansions up to the second order, we obtain
\begin{align*}
|\E f(S)-f(Z)|&\le\|f'\|_\infty\E|S-Z|\\
\mbox{and}\quad|\E f(Y)-f(\sigma N)|&\le\frac{\|f''\|_\infty^2}{2}\E|Y-\sigma N|^2.
\end{align*}
We note that $Z-S= (V_1+\cdots+V_k)/\sqrt{n}$ is a sum of $X_i$'s
for which the number of terms is $\le (k+1)q+p$. Then \refeqn{covepsilon} and \refeqn{covkappa}, under conditions (\ref{varepsilon}) or
(\ref{varkappa}) respectively,
entail:
$$(\E|Z-S|)^2\le\E|Z-S|^2\preceq \big((k+1)q+p\big)/n. $$
Now $Y\sim\displaystyle\sqrt\frac{kp}{n}\sigma_p N$, thus
\begin{equation*}
\E|Y-\sigma
N|^2\le\left|\frac{kp}{n}-1\right|\sigma_p^2+|\sigma_p^2-\sigma^2|.
\end{equation*}
Remarking that $\left|kp/n-1\right|^2\preceq \big((k+1)q+p\big)/n$, it remains to bound the quantity
\begin{eqnarray*}
|\sigma_p^2-\sigma^2|&\le&\sum_{|i|<p}\frac{|i|}p|\E
X_0X_i|+\sum_{|i|>p}|\E X_0X_i|.
\end{eqnarray*}
Let $a_i=|\E X_0X_i|$, under conditions (\ref{varepsilon}) or
(\ref{varkappa}) (respectively), the series $\displaystyle \sum_{i=0}^\infty a_i$
converge thus $\displaystyle s_j=\sum_{i=j}^\infty a_i\cvg\limits_{j\to\infty}0$ and
\begin{eqnarray*}
|\sigma_p^2-\sigma^2| \le2\sum_{i=0}^{p-1}\frac{i}p\cdot a_i+2s_p
\le\frac2p\sum_{i=0}^{p-1}s_i+2s_p.
\end{eqnarray*}
Cesaro lemma entails that term $|\sigma_p^2-\sigma^2|$ converges to
0.\\
Hence $|\E f(S)-f(Z)|+|\E f(Y)-f(\sigma N)|$ tends to 0 as
$n\uparrow \infty $.

\noindent To determine the convergence rate, we assume that $a_i={\cal
O}(i^{-\alpha})$ for some $\alpha>1$;  then
$$
|\sigma_p^2-\sigma^2|\preceq p^{1-\alpha}.
$$
Remarking that $a_i=\E X_0X_i=\cov(X_0,X_i)$, we then use equations
(\ref{covkappa}) and (\ref{covepsilon}) and we find $\alpha=\kappa$
or $\alpha=\lambda(m-2)/(m-1)$ depending of the weak-dependence
setting. With $p=n^a$, $q=n^b$ for 2 constants $a$ and $b$ and from
the relation $\|f^{(j)}\|_\infty\le |t|^j$, those bounds become, up
to a constant
\begin{description}
\item
$|t|\left(n^{(b-a)/2}+n^{(a-1)/2}\right)+t^2\left(n^{b-a}+n^{a(1-\kappa)}\right),$
in the $\kappa$-weak dependence setting,

\item $|t|\left(n^{(b-a)/2}+n^{(a-1)/2}\right)+t^2\left(n^{b-a}+n^{a(1-\lambda(m-2)/(m-1))}\right),$  for $\lambda$-weak
dependence.
\end{description}
\subsubsection{Main terms}
It remains to control the second and the third terms of the sum.
They are bounded as usual by
\begin{eqnarray*}
|\E f(Z)-f(Z^*)|\le\sum_{j=1}^k|\E\Delta_j|,\quad |\E
f(Z^*)-f(Y)|\le\sum_{j=1}^k|\E\Delta '_j|,
\end{eqnarray*}
where $\Delta_j=f(W_j+x_j)-f(W_j+x^*_j)$, for $j=1,\ldots,k$  with
$x_j=\frac1{\sqrt{n}}{U_j},$ $x^*_j=\frac1{\sqrt{n}}{U^*_j}$,
$W_j=w_j+\sum_{i>j}x^*_i$, $w_j=\sum_{i<j}x_i$ and $\Delta
'_j=f(W'_j+x^*_j)-f(W'_j+x'_j)$, for $j=1,\ldots,k$ with
$x'_j=\frac1{\sqrt{n}}{U'_j}$,
$W'_j=\sum_{i<j}x^*_i+\sum_{i>j}x'_i.$
\\
Exploiting the special form of $f$ and the independence properties
of the variables $U^*_i$ and $U'_i$, we can write
\begin{eqnarray*}
\E\Delta_j&=&\left(\E f(w_j)f(x_j)-\E f(w_j)\E f(x_j^*)\right)\E
f\left(\sum_{i>j}x^*_i\right),\\
\E\Delta'_j&=&\left(\E f(x_j^*)-\E f(x'_j)\right)\E
f\left(W'_j\right).
\end{eqnarray*}
We then control the two terms $\displaystyle \E
f\left(\sum_{i>j}x^*_i\right)$ and $\displaystyle \E f\left(W'_j\right)$ by the
fact that $\|f\|_\infty\le 1$ and we use the coupling to introduce a
covariance term
\begin{eqnarray*}
|\E\Delta_j|&\le&\left|\cov\left(f\left(\sum_{i<j}x_i\right),f(x_j)\right)\right|,\\
|\E\Delta'_j|&\le&\left|\E f(x_j^*)-\E f(x'_j)\right|.
\end{eqnarray*}

\begin{itemize}
\item For $\Delta_j$, we use weak dependence.\\
To do so, write $|\E\Delta_j|=|\cov[F(X_m,m\in B_i, i<j),G( X_m,m\in
B_j)]|$, with $\displaystyle
F(z_1,\ldots,z_{kp})=f\Big(\sum_{i<j}u_i/\sqrt{n}\Big)$ where $u_i=\sum_{\ell\in B_i}z_\ell$. We verify that $\|F\|_\infty\le 1$ and we control $\Lip F$:
\begin{multline*}
\left|f\left(\frac1{\sqrt{n}}\sum_{i<j}\sum_{\ell\in
B_i}z_\ell\right)-f\left(\frac1{\sqrt{n}}\sum_{i<j}\sum_{\ell\in
B_i}z'_\ell\right)\right|\\
\le
\left|1-\exp{it\left(\frac1{\sqrt{n}}\sum_{i<j}\sum_{\ell\in B_i}(z_\ell-z'_\ell )\right)}\right|\\
\le\frac{|t|}{\sqrt{n}}\sum_{\ell=1}^{kp}|z_\ell-z'_\ell|.
\end{multline*}
Similarly, for $\ds G(z_1,\ldots,z_p)=f\left(\sum_{i=1}^pz_i/\sqrt{n}\right)$, we have $\|G\|_\infty=1$ and $\Lip G\preceq |t|/\sqrt{n}$. We then distinguish the two cases of weak dependence, remarking the gap between the left and the right terms in the covariance is at least
$q$.
\begin{itemize}
\item In the $\kappa$-weak dependent setting:
$|\E\Delta_j|\preceq kp\cdot p\cdot \frac{|t|}{\sqrt{n}}\cdot
\frac{|t|}{\sqrt{n}}\cdot \kappa(q).$
\item Under the $\lambda$ dependence condition:\\
$|\E\Delta_j|\preceq \left(kp\cdot p\cdot \frac{|t|}{\sqrt{n}}\cdot
\frac{|t|}{\sqrt{n}}
+kp\cdot\frac{|t|}{\sqrt{n}}+p\cdot\frac{|t|}{\sqrt{n}}\right)\cdot\lambda(q).$
\end{itemize}
Note that these bounds do not depend on $j$:
\begin{align*}
|\E f(Z)-f(Z^*)|
&\preceq kp\cdot t^2\cdot\kappa(q),&\mbox{ under } \kappa,\\
&\preceq kp\cdot(t^2+|t|\sqrt{k/p})\cdot\lambda(q),& \mbox{
under } \lambda.
\end{align*}
Knowing that $p=n^a$, $q=n^b$, $\kappa(r)={\cal
O}\left(r^{-\kappa}\right)$ or $\lambda(r)={\cal
O}\left(r^{-\lambda}\right)$, these convergence rates respectively become $n^{1-\kappa b}$ or
$n^{1+(1/2-a)_+-\lambda b}$ in the $\kappa$ or the
$\lambda$ dependence context.
\item For $\Delta'_j$, Taylor expansions up to order 2 or 3 respectively give:
\begin{eqnarray*}
|f(x^*_j)-f(x'_j)|&\le&|x^*_j-x'_j|\|f'\|_\infty+\frac12(x^*_j-x'_j)^2\|f''\|_\infty+r_j\\
r_j&\le &\frac12\|f''\|_\infty (x^*_j-x'_j)^2, \mbox{ or}\\
&\le& \frac16\|f'''\|_\infty|x^*_j-x'_j|^3,
\end{eqnarray*}
For an arbitrary $\delta\in[0,1]$, we have:
\begin{eqnarray*}
\E r_j&\preceq&\E(t^2(|x^*_j|^2+|x'_j|^2)\wedge |t|^3(|x^*_j|^3+|x'_j|^3))\\
&\preceq&
\E(t^2|x^*_j|^2\wedge |t|^3|x^*_j|^3)+\E(t^2|x'_j|^2\wedge |t|^3|x'_j|^3)\\
&\preceq&
|t|^{2+\delta}\left(\E|x^*_j|^{2+\delta}+\E|x'_j|^{2+\delta}\right).
\end{eqnarray*}
By the stationarity of the sequence $(X_i)_{i\in\Z}$, we
obtain
\begin{eqnarray*}
|\E \Delta '_j|
&\preceq&|t|^{2+\delta}n^{-1-\frac\delta2}\left(\E|S_p|^{2+\delta}\vee
p^{1+\frac\delta2}\right).
\end{eqnarray*}
Lemma \ref{momlambda} allows us to fin a bound for
$\E|S_p|^{2+\delta}$. If $\kappa>2+\frac{1}{\zeta}$, or
$\lambda>4+\frac{2}{\zeta},$ where $\kappa(r)={\cal
O}\left(r^{-\kappa}\right)$ or $\lambda(r)={\cal
O}\left(r^{-\lambda}\right)$ then there exist $\delta\in]0,
\zeta\wedge1[$ and $C>0$ such that
$$
\E|S_p|^{2+\delta}\le Cp^{1+\delta/2}.
$$
We then obtain
$$
|\E f(Z^*)-f(Y)|\preceq
|t|^{2+\delta}k(p/n)^{1+\delta/2}.
$$
Because $p=n^a$, this bound is of order $n^{(a-1)\delta/2}$ in both $\kappa$ and $\lambda$-weak dependence settings.
\end{itemize}
We now collect the previous bounds to conclude that a
multidimensional CLT holds under assumptions of both Theorems
\ref{th1} and \ref{th2}. Tightness follows from the Kolmogorov-Chentsov
criterion (see \cite{B}) and Lemma \ref{momlambda}; thus both Theorems
\ref{th1} and \ref{th2} follow from repeated application of the
previous CLT.\ $\square$
\subsection{Rates of convergence}
\label{rate}
Rates of convergence are  now presented in two propositions of independent interest.
We compute explicit bounds for both the difference of characteristic functions and the  Berry-Ess\'een inequalities.
\begin{proposition}
\label{prop2}
Let $(X_t)_{t\in\Z}$ be a weakly dependent stationary process satisfying \refeqn{mo2+} with $m=2+\zeta$ then the difference between the
characteristic functions is bounded by
$$
\left|\E \left(e^{itS_n/\sqrt n}-e^{it\sigma N}\right)\right|=o(n^{-c}),
$$
for some $c<c^*$  and all $t\in\R$  where
$c^\ast$ depends of the weak dependent coefficients
\begin{itemize}
\item under $\kappa$-weak dependence, if $\kappa(r)={\cal O}(r^{-\kappa})$ for $\kappa>2+\frac{1}{\zeta}$,
then we have $\displaystyle c^*=\frac{(\kappa-1)A}{A+2\kappa(1+A)}$ where
$$
\displaystyle
A=\frac{\sqrt{(2\kappa-3-\zeta)^2+4(\kappa\zeta-2\zeta-1)}+\zeta+3-2\kappa}{2}\wedge
1.
$$
\item under $\lambda$-weak dependence, if $\lambda(r)={\cal O}(r^{-\lambda})$ for
$\lambda>4+\frac{2}{\zeta}$, then we obtain $\displaystyle
c^*=\frac{(\lambda+1)B}{2+B+2\lambda(1+B)}$ where
$$
\displaystyle
B=\frac{\sqrt{(2\lambda-6-\zeta)^2+4(\lambda\zeta-4\zeta-2)}+\zeta+6-2\lambda}{2}\wedge
1,
$$
\end{itemize}
\end{proposition}
We  use the following Esséen inequality in Proposition~\ref{prop3}
\begin{theorem}[Theorem 5.1 p.142 of \cite{Petrov}]
\label{theo:petrov}
Let $X$ and $Y$ be 2 random variables  and assume that $Y$ is Gaussian. Let $F$ and $G$
be their distribution functions with corresponding characteristic
functions $f$ and $g$. Then, for every $T>0$, we have for suitable
constants $b$ and $c$
\begin{equation}%
\label{petrov}
\sup_{x\in\R}|F(x)-G(x)|\le
b\int_{-T}^T\left|\frac{f(t)-g(t)}{t}\right|dt+ \frac{c}{T}.
\end{equation}
\end{theorem}
\begin{proposition}[A rate in the Berry Essen bounds]
\label{prop3}
Let $(X_t)_{t\in\Z}$ be a real stationary process
satisfying Proposition \ref{prop2} assumptions. We obtain
$$
\sup_x|F_n(x)-\Phi(x)|=o\left(n^{-c}\right)
$$
with $c<c'$ where $c'=c^\ast/(3+A)$ or
$c'=c^\ast/(3+B)$, respectively, in $\kappa$ or $\lambda$-weak dependence contexts
($A$, $B$ and $c^\ast$ are defined in Proposition \ref{prop2}).
\end{proposition}
\begin{proofof}{Proposition \ref{prop2}}
 In the previous
section, the different terms have already been bounded as follows:
\begin{itemize}
 \item In the
$\kappa$-weak dependence case, the exponents of $n$ in the bounds
obtained in section \ref{sectiontheo12} are
\begin{itemize}
\item for the auxiliary terms: $(b-a)/2$, $(a-1)/2$ and $a(1-\kappa)$, \item for the main
terms: $1-\kappa b$ and $(a-1)\delta/2$.
\end{itemize}
Because $\delta<1$ and $b<a<1$, we remark that
$(a-1)\delta/2>(a-1)/2$ and $1-\kappa b>a(1-\kappa)$. The only rate
of the auxiliary term it remains to consider is $(b-a)/2$ and we
obtain
\begin{eqnarray*}
a^*=\frac{2+\delta+2\kappa\delta}{\delta+2\kappa(1+\delta)}\in\left]b^\ast,\frac{\delta}{1+\delta}\right[,\
\ b^*=\frac{2+a^\ast}{1+2\kappa}\in\left]0,a^*\right[.
\end{eqnarray*}
We conclude with standard calculations and with the help of the
inequality $\delta<A$ (see \refeqn{A}).

\item We have the equivalent in the $\lambda$-weak dependence case
\begin{itemize}
\item for the auxiliary terms: $(b-a)/2$, $(a-1)/2$ and $a(1-\lambda)$, \item for the main
terms: $1+(1/2-a)_+-\lambda b$ and $(a-1)\delta/2$.
\end{itemize}
Only three rates give the asymptotic: $(a-1)\delta/2$,
$1+(1/2-a)_+-\lambda b$ and $(b-a)/2$. In the previous case, the
optimal choice of $a^*$ was smaller than $1/2$. Then we have to
consider here the rate $3/2-a-\lambda b$ and not $1-\lambda b$. Thus
\begin{align*}
a^*&=\frac{3+\delta+2\lambda\delta}{2+\delta+2\lambda(1+\delta)}\in\introo{b^\ast}{\frac{\delta}{1+\delta}},\\
b^*&=\frac{3+2\delta}{2+\delta+2\lambda(1+\delta)}\in\introo{0}{a^*}
\end{align*}
\end{itemize}
Finally, we obtain a rate of $n^{-c^*}$ using the inequality
\refeqn{B}.
\end{proofof}
\begin{proofof}{Proposition \ref{prop3}}
Let choose $a^\ast$ and $b^\ast$ as in the proof of proposition \ref{prop2}. We now need to make precise the
impact of $t$ on the different term of the bound of
the $\L^1$ distance between
the characteristic functions of $S$ and $\sigma N$. Up to a
constant independent of $t$, the Kolmogorov distance is bounded by $\left(|t|+t^2+|t|^{2+C}\right)n^{-c^\ast}$. Here $C=A$ or
$B$ in the two contexts of dependence. Using Theorem \ref{theo:petrov} for a well chosen value of $T$, we
obtain the result of proposition \ref{prop3}.
\end{proofof}
\subsection{Proof of Lemma \ref{lemlipsch}}\label{prooflipsch}
The case of Lipschitz functions of dependent inputs is divided in
two sections devoted respectively to the definition of such models
and to their weak dependence properties.
\subsubsection{Existence}
Let $Y^{(s)}=(Y_{-i}\1_{|i|<s})_{i\in\Z}$,
$Y^{(s)}_+=(Y_{-i}\1_{-s<i\le s})_{i\in\Z}$ for $s\in\Z$ and
$H\left(Y^{(\infty)}\right)=\lim_{s\to\infty}H\left(Y^{(s)}\right)$.
In order to prove the existence of the Bernoulli shift with dependent
inputs, we show that $X_0$ is the sum of a normally convergent
series in $\L^m$; formally
\begin{multline*}
X_0=H\big(Y^{(\infty)}\big)=H(0)+\left(H\big(Y^{(1)}\big)-H(0)\right)\\
+\sum_{s=1}^\infty
H\big(Y^{(s+1)}\big)-H\big(Y^{(s)}_+\big)+\left(H\big(Y^{(s)}_+\big)-H\big(Y^{(s)}\big)\right).
\end{multline*}
>From \refeqn{as} we obtain
\begin{eqnarray*}
\left\|H\left(Y^{(1)}\right)-H(0)\right\|_m&\le& b_0\|Y_0\|_m,\\
\left\|H\left(Y^{(s+1)}\right)-H\left(Y^{(s)}_+\right)\right\|_m&\le& b_{-s}\|Y_{-s}\|_m,\\
\left\|H\left(Y^{(s)}_+\right)-H\left(Y^{(s)}\right)\right\|_m&\le&
b_{s}\|Y_{s}\|_m.
\end{eqnarray*}
By $(Y_t)_{t\in\Z}$'s stationarity we get
\begin{multline}\label{as1}
\|X_0\|_m\le
\left\|H\left(Y^{(1)}\right)-H(0)\right\|_m+\sum_{s=1}^\infty
\left\|H\left(Y^{(s+1)}\right)-H\left(Y^{(s)}_+\right)\right\|_m\\
+\left\|H\left(Y^{(s)}_+\right)-H\left(Y^{(s)}\right)\right\|_m\le \sum_{i\in\Z} b_i\|Y_0\|_m
\end{multline}
Analogously, the process $X_t=H(Y_{t-i},i\in\Z)$ is well defined as
the sum of a normally convergent series in $\L^m$. The stationarity
of $(X_t)_{t\in\Z}$ holds from that of the input process
$(Y_t)_{t\in\Z}$.
\subsubsection{Weak dependence properties}
Let $X_n^{(r)}=H(Y^{(r)})$ and $X_{\bf s}=(X_{s_1},\dots,X_{s_u})$, $X_{\bf t}=(X_{t_1},\dots,X_{t_v})$ for
any $k\ge 0$ and any $(u+v)$-tuple such that $s_1<\cdots< s_u\le
s_u+k\le t_1<\cdots<t_v$. Then we have for all $f,g$ satisfying
$\|f\|_\infty,\|g\|_\infty\le1$ and $\Lip f+\Lip g<\infty$
\begin{align}
|\cov(f(X_{\bf s}),g(X_{\bf t}))|
&\le|\cov(f(X_{\bf s})-f(X_{\bf s}^{(r)}),g(X_{\bf t}))|
\label{term1}\\
&\qquad+~|\cov(f(X_{\bf s}^{(r)}),g(X_{\bf t})-g(X_{\bf t}^{(r)}))|
\label{term2}\\
&\qquad+~|\cov(f(X_{\bf s}^{(r)}),g( X_{\bf t}^{(r)}))|.
\label{term3}
\end{align}
Using the fact that $\|g\|_\infty\le1$, we bound the term in \refeqn{term1}
$$2\Lip
f\cdot \E\left|\sum_{i=1}^u(X_{s_i}-
X_{s_i}^{(r)})\right|\le2u\Lip f \max_{1\le i\le
u}\E\left|X_{s_i}-X_{s_i}^{(r)}\right|.
$$
Applying inequality (\ref{as1}) in the case where $m=1$,
we obtain
$\E\left|X_{s_i}-X_{s_i}^{(r)}\right|\le\sum_{i\ge r}b_i\|Y_0\|_1$.
The second term (\ref{term2}) is bounded in a similar way.

\noindent The last term (\ref{term3}) can be written as
$$\left|\cov(F^{(r)}(Y_{s_i+j},1\le i\le u,|j|\le
r),G^{(r)}(Y_{t_i+j},1\le i\le v,|j|\le r)\right|,
$$
where $F^{(r)}:\R^{u(2r+1)}\to\R$ and $G^{(r)}:\R^{v(2r+1)}\to\R$.
Under the assumption $r\le [k/2]$,
we use the $\epsilon=\eta$ or $\lambda$-weak dependence of $Y$ in order to bound this covariance term by $\psi(\Lip F^{(r)},\Lip
G^{(r)},u(2r+1),v(2r+1))\epsilon_{k-2r}$,
with respectively $\psi(u,v,a,b)=ua+vb$ or $\psi(u,v,a,b)=uvab+ua+vb$.
We compute
\begin{multline*}
\Lip F^{(r)}=\\
\sup \frac{|f(H(x_{s_i+l},1\le i\le u,|l|\le
r)-f(H(y_{s_i+l},1\le i\le u,|l|\le r)|}{\sum_{i=1}^u\sum_{-r\le
l\le r}|x_{s_i+l}-y_{s_i+l}|},
\end{multline*}
where the $\sup$ extends to $x\ne y$ where $x,y\in\R^{u(2r+1)}$. Notice
now that if $x, y$ are sequences with $x_i=y_i=0$ if $|i|\ge r$ then
repeated applications of the condition (\ref{as}) yields
\begin{equation}%
\label{as21}
|H(x)-H(y)|\le \sum_{|i|\le r}b_i|x_i-y_i|\le L \sum_{|i|\le
r}|x_i-y_i|\end{equation} where $L=\sum_{i\in\Z}b_i$. Repeating
inequality \refeqn{as21}, we obtain
$$
|F^{(r)}(x)- F^{(r)}(y)| \le \Lip f\, L \sum_{i=1}^u\sum_{-r\le l\le
r}|x_{s_i+l}-y_{s_i+l}|
$$
and we get $\Lip F^{(r)}\le \Lip f \cdot L$. Similarly $\Lip
G^{(r)}\le \Lip g \, L$.\\ Under $\eta$-weak dependent inputs, we
bound the covariance
\begin{multline*}
|\cov(f(X_{\bf s}),g(X_{\bf t}))|\\\le (u\Lip f+v\Lip g) \times\left[
2\sum_{|i|\ge r}b_i\|Y_0\|_1+(2r+1)L\eta_{Y}(k-2r)\right].
\end{multline*}
Under $\lambda$-weak dependent inputs
\begin{multline*}
|\cov(f(X_{\bf s}),g(X_{\bf t}))|
\le (u\Lip f+v\Lip g+uv\Lip f\Lip g)\times\\
\times\Big(\Big\{2\sum_{|i|\ge
r}b_i\|Y_0\|_1+(2r+1)L\,\lambda_{Y}(k-2r)\Big\}\vee \\
(2r+1)^2L^2\lambda_{Y}(k-2r)\Big). \ \square\end{multline*}
\subsection{Proof of Lemma \ref{lambd}}\label{prooflemma}
\subsubsection{Existence}
We decompose $X_0$ as above in the case $\ell=0$. Here, we bound each terms by
\begin{eqnarray*}
|H(Y^{(1)})-H(0)|&\le& b_0|Y_0|\\
|H(Y^{(s+1)})-H(Y^{(s)}_+)|&\le& b_{-s}(\|Y^{(s)}_+\|_\infty^l\vee 1)|Y_{-s}|\\
|H(Y^{(s)}_+)-H(Y^{(s)}|&\le& b_{s}(\|Y^{(s)}\|_\infty^l\vee
1)|Y_{s}|
\end{eqnarray*}
Using H\"older inequality yields
\begin{multline*}
\E\left|H(Y^{(1)})-H(0)\right|+\sum_{s=1}^\infty
\E\left|H(Y^{(s+1)})-H(Y^{(s)}_+)\right|\\
+\E\left|H(Y^{(s)}_+)-H(Y^{(s)})\right|
\le \sum_{i\in\Z}2|i| b_i(\|Y_0\|_1+\|Y_0\|_{l+1}^{l+1})
\end{multline*}
Hence assumptions $\ell+1\le {m'} $ and $\sum_{i\in\Z}|i| b_i<\infty$
together imply that the variable $H(Y)$ is well defined in $\L^1$.
In the same manner, the process $X_n=H(Y_{n-i},i\in\Z)$ is well
defined.
The proof extends in $\L^m$ if $m\ge 1$ is such that $(\ell+1)m\le {m'}$.

\subsubsection{Weak dependence properties}
Here, we exhibit some Lipschitz functions and
we then truncate inputs. We write $\Ybar=Y\vee(-T)\wedge T$ for a truncation
$T$ set below. Denote $X_n^{(r)}=H(Y^{(r)})$ and $\Xbar_n^{(r)}=H(\Ybar^{(r)})$. Furthermore, for any $k\ge 0$ and
any $(u+v)$-tuple such that $s_1<\cdots< s_u\le s_u+k\le
t_1<\cdots<t_v$, we set $X_{\bf s}=(X_{s_1},\dots,X_{s_u})$, $X_{\bf
t}=(X_{t_1},\dots,X_{t_v})$ and $\Xbar_{\bf
s}^{(r)}=(\Xbar_{s_1}^{(r)},\dots,\Xbar_{s_u}^{(r)})$,
$\Xbar_{\bf t}^{(r)}=(\Xbar_{t_1}^{(r)},\dots,\Xbar_{t_v}^{(r)})$. Then we have for all $f,g$ satisfying
$\|f\|_\infty,\|g\|_\infty\le1$ and $\Lip f+\Lip g<\infty$
\begin{align}
|\cov(f(X_{\bf s}),g(X_{\bf t}))|
&\le|\cov(f(X_{\bf s})-f(\Xbar_{\bf s}^{(r)}),g(X_{\bf t}))|
\label{term12}\\
&\qquad+~|\cov(f(\Xbar_{\bf s}^{(r)}),g(X_{\bf t})-g(\Xbar_{\bf t}^{(r)}))|
\label{term22}\\
&\qquad+~|\cov(f(\Xbar_{\bf s}^{(r)}),g(\Xbar_{\bf t}^{(r)}))|.
\label{term32}
\end{align}
Using the fact that $\|g\|_\infty\le1$, the term (\ref{term12}) is bounded by
$$
2u\Lip f \left(\max_{1\le i\le
u}\E\left|X_{s_i}-X_{s_i}^{(r)}\right|+\max_{1\le i\le
u}\E\left|X_{s_i}^{(r)}-\Xbar_{s_i}^{(r)}\right|\right).
$$
With the same arguments used in the proof of the existence of $H(Y^{(\infty)})$,
the first term in the right-hand side of the inequality is bounded
by
$$
\sum_{i\ge s}2|i| b_i(\|Y_0\|_1+\|Y_0\|_{l+1}^{l+1}).
$$
Notice now that if $x, y$ are sequences with $x_i=y_i=0$ if $|i|\ge r$ then an infinitely repeated application of the previous inequality (\ref{as}) yields
\begin{equation}%
\label{as2}
|H(x)-H(y)|\le L(\|x\|_\infty^l\vee\|y\|_\infty^l\vee
1)\|x-y\|
\end{equation}
where $L=\sum_{i\in\Z}b_i<\infty$ because
$\sum_{i\in\Z}|i| b_i<\infty$. The second term is bounded by
using \refeqn{as2}
\begin{eqnarray*}\E\left|X_{s_i}^{(r)}-\Xbar_{s_i}^{(r)}\right|&=& \E\left|H\left(Y^{(r)}\right)-H\left(\Ybar^{(r)}\right)\right|\\
&\le&L\E\left(\left(\max_{-r\le i\le r}|Y_i|\right)^l\sum_{-r\le
j\le
r}|Y_j|\1_{Y_j\ge T}\right)\\
&\le&L(2r+1)^2\E\left(\max_{-r\le i,j\le
r}|Y_i|^l|Y_j|\1_{|Y_j|\ge T}\right)\\
&\le&L(2r+1)^2\|Y_0\|_{m'}^{m'}T^{\ell+1-{m'}}
\end{eqnarray*}
The second term (\ref{term22}) of the sum is analogously bounded. The
last term (\ref{term32}) can be written as
$$
\left|\cov\left(\Fbar^{(r)}(Y_{s_i+j},1\le i\le u,|j|\le
r),\Gbar^{(r)}(Y_{t_i+j},1\le i\le v,|j|\le r)\right)\right|,
$$
where
$\Fbar^{(r)}:\R^{u(2r+1)}\to\R$ and $\Gbar^{(r)}:\R^{u(2r+1)}\to\R$. Under the assumption $r\le [k/2]$, we
use the $\epsilon=\eta$ or $\lambda$-weak dependence of $Y$ in order
to bound this covariance term by $\psi\left(\Lip \Fbar^{(r)},\Lip
\Gbar^{(r)},u(2r+1),v(2r+1)\right)\epsilon_{k-2r}$,  we set
respectively $\psi(u,v,a,b)=uvab$ or $\psi(u,v,a,b)=uvab+ua+vb$.
\begin{multline*}
\Lip \Fbar^{(r)}=\\
\sup \frac{|f(H(\xbar_{s_i+l},1\le i\le
u,|l|\le r)-f(H(\ybar_{s_i+l},1\le i\le u,|l|\le
r)|}{\sum_{j=1}^u\|x_j-y_j\|},
\end{multline*}
where the $\sup$ extends to $(x_1,\ldots, x_u)\ne (y_1,\ldots, y_u)$
where $x_i,y_i\in\R^{2r+1}$. Using \refeqn{as2}
\begin{eqnarray*}
|\Fbar^{(r)}(x)-\Fbar^{(r)}(y)| &\le& \Lip f L
\sum_{i=1}^u\left(\|\xbar_{s_i}\|_\infty\vee\|\ybar_{s_i}
\|_\infty\vee1\right)^l\|\xbar_{s_i}-\ybar_{s_i}\|\\
&\le& \Lip f L T^l\sum_{i=1}^u\sum_{-r\le l\le
r}|x_{s_i+l}-y_{s_i+l}|.
\end{eqnarray*}
We thus obtain  $\Lip F^{(r)}\le \Lip f \cdot L\cdot T^l $.
Similarly $\Lip G^{(r)}\le \Lip g \cdot L\cdot T^l $.\\ Under
$\eta$-weak dependent inputs, we bound the covariance
\begin{multline*}
|\cov(f(X_{\bf s}),g(X_{\bf t}))|\le (u\Lip f+v\Lip g)\Big\{ 4\sum_{|i|\ge r}|i|
b_i(\|Y_0\|_1+\|Y_0\|_{l+1}^{l+1})\\
+(2r+1)L\,\big((2r+1)2\|Y_0\|_{m'}^{m'}T^{l+1-{m'}}
+T^l\eta_{Y}(k-2r)\big)\Big\}
\end{multline*} We then fix the truncation $\displaystyle
T^{{m'}-1}=\frac{2(2r+1)\|Y_0\|_{m'}^{m'}}{\eta_{Y}(k-2r)}$ to
obtain the result of Lemma \ref{lambd} in the $\eta$-weak
dependent case.\\ Under $\lambda$-weak dependent inputs
\begin{multline*}
|\cov(f(X_{\bf s}),g(X_{\bf t}))|\le (u\Lip f+v\Lip g+uv\Lip f\Lip g)\times\\
\times\Big(\Big\{ 4\sum_{|i|\ge r}|i|
b_i(\|Y_0\|_1+\|Y_0\|_{l+1}^{l+1})\\
+(2r+1)L\left(2(2r+1)T^{l+1-{m'}}\|Y_0\|_{m'}^{m'}+T^l\lambda_{Y}(k-2r)\right)\Big\}
\\
\qquad\vee \ \left\{(2r+1)^2L^2T^{2l}\lambda_{Y}(k-2r)\right\}\Big)
\end{multline*}
We then set a truncation such that $\displaystyle
T^{l+{m'}-1}=\frac{2\|Y_0\|_{m'}^{m'}}{L\lambda_{Y}(k-2r)}$ to
obtain the result of Lemma \ref{lambd} in the $\eta$-weak
dependent case.\ $\square$

\acknowledgements We wish to thank  anonymous referees for his
comments; the first one pointed out the reference of \cite{heyde};
we
 we then exhibited examples for which this method based on
martingale arguments does not seem to apply. A  second referee  helped us to improve on the readability of the manuscript. We are deeply grateful to Alain Latour who made a critical review of the drafts and with whom we have worked on the final version of this paper.

\bibliography{TLCdoukwint}
\end{document}